\documentclass{amsart}

\usepackage{amsmath}
\usepackage{amsthm}
\usepackage{amssymb}
\usepackage{xspace}
\usepackage{color}
\usepackage{mathtools}

\theoremstyle{plain}

\newtheorem{theorem}{Theorem}
\newtheorem{lemma}[theorem]{Lemma}

\theoremstyle{definition}

\newtheorem{example}[theorem]{Example}

\newcommand{\comment}[1]{}
\newcommand{\rea}{\ensuremath{\mathbf{R}}\xspace}
\newcommand{\probs}{\ensuremath{\mathcal{P}}\xspace}
\newcommand{\intd}[1]{\,\mathrm{d}#1 \,}
\newcommand{\closure}[1]{\ensuremath{\overline{#1}}\xspace}
\newcommand{\dimh}{\ensuremath{\dim_{\mathrm{H}}}\xspace}
\newcommand{\dimf}{\ensuremath{\dim_{\mathrm{F}}}\xspace}
\newcommand{\dimfc}{\ensuremath{\dim_{\mathrm{FC}}}\xspace}
\newcommand{\leb}{\ensuremath{\lambda}\xspace}

\newcommand\restr[2]{{							
  \left.\kern-\nulldelimiterspace				
  #1											
  \right|_{#2}									
  }}

\DeclareMathOperator{\re}{Re}
\DeclareMathOperator{\sinc}{sinc}
\DeclarePairedDelimiter{\ceil}{\lceil}{\rceil}

\begin{document}
\title{The Fourier dimension is not finitely stable}

\author{Fredrik Ekstr\"om}
\address{Centre for Mathematical Sciences\\ Lund University\\ Box 118\\ 22 100 Lund\\ SWEDEN}
\email{fredrike@maths.lth.se}

\begin{abstract}
The Fourier dimension is not in general stable under finite unions of sets. Moreover,
the stability of the Fourier dimension on particular pairs of sets is independent
from the stability of the compact Fourier dimension.
\end{abstract}

\maketitle

\section{Introduction}
\noindent
The present note gives an example to show that the Fourier dimension is not
stable under finite unions of sets. This improves one of the results in
\cite{EPS}, where it was shown that the Fourier dimension is not countably
stable. For a more detailed introduction to the Fourier dimension than is
given here, see for example \cite{EPS} and references therein.

The \emph{Fourier transform} of a finite Borel measure $\mu$ on $\rea^d$ is
defined as
	$$
	\widehat\mu(\xi) = \int e^{-2\pi i \xi \cdot x} \intd{\mu}(x),
	$$
where $\xi \in \rea^d$ and $\cdot$ denotes the Euclidean inner product.
The \emph{Fourier dimension} of $\mu$ measures how quickly $\widehat\mu$
decays at infinity, and is defined by
	$$
	\dimf \mu = \sup\left\{ s \in [0, d]; \, \exists C \in \rea \text{ such that }
	|\widehat\mu(\xi)| \leq C |\xi|^{-s/2} \text{ for all } \xi \right\}.
	$$
If $A$ is a Borel subset of $\rea^d$ then the Fourier dimension of $A$ is defined
to be
	$$
	\dimf A = \sup\left\{\dimf \mu; \, \mu \in \mathcal{P}(A) \right\},
	$$
where $\probs(A)$ denotes the set of Borel probability measures on $\rea^d$
that give full measure to $A$. It can be shown (see \cite[Lemma~12.12]{mattila})
that if $\mu$ has compact support and $0 < \dimf \mu < d$ then
	$$
	I_s(\mu) := \iint |x - y|^{-s} \intd{\mu}(x) \intd{\mu}(y)
	$$
is finite for all $s < \dimf \mu$, and from this it follows that
$\dimf A \leq \dimh A$ for any Borel set $A$.

As an alternative way to define the Fourier dimension, one could require that the
measure in the definition should give full measure to some compact subset of $A$
and not just to $A$ itself. This \emph{compact Fourier dimension} is thus defined
as
	$$
	\dimfc A = \sup \left\{ \dimf \mu; \, \mu \in \probs(K), \, K \subset A \text{ is compact} \right\}.
	$$
One of the examples in \cite{EPS} shows that the compact Fourier dimension is
not finitely stable, but that example is not a counterexample to finite stability
of the Fourier dimension.
Example~\ref{theexample} below produces the opposite situation, namely, sets
$A'$ and $B'$ such that
	\begin{align*}
	\dimf (A' \cup B') &> \max(\dimf A', \, \dimf B'),\\
	\dimfc (A' \cup B') &= \max(\dimfc A', \, \dimfc B').
	\end{align*}

\section{The example}
\noindent
The following lemma is used in the example. It previously appeared in \cite{EPS},
but is included here for completeness.

\begin{lemma}	\label{effinfsuplemma}
For any $\varepsilon \in (0, 1]$,
	$$
	\inf_{\phantom{!}\mu\phantom{!}} \sup_{j \geq 1} \left| \widehat \mu(j) \right|
	\geq \frac{\pi\varepsilon}{8 + 2\pi\varepsilon}
	\qquad\left( \, \geq \frac{\varepsilon}{5} \, \right),
	$$
where the infimum is over all $\mu \in \probs([\varepsilon, 1])$ and
the supremum is over all positive integers $j$.

\begin{proof}
Fix $\varepsilon > 0$ and take any $\mu \in \probs([\varepsilon, 1])$. If $\varphi$
is a real-valued continuous function supported on $[0, \varepsilon]$ such
that
	$$
	\int \varphi(x) \intd{x} = 1 \qquad \text{ and } \qquad
	\sum_{k = -\infty}^\infty |\widehat\varphi(k)| < \infty
	$$
%
then
	$$
	0 = \mu(\varphi) = \sum_{k = -\infty}^\infty \widehat\varphi(k) \overline{\widehat\mu(k)}
	= 1 + 2 \re\left( \sum_{k = 1}^\infty \widehat\varphi(k) \overline{\widehat\mu(k)} \right),
	$$
and thus
	$$
	\frac{1}{2} \leq \sum_{k = 1}^\infty |\widehat\varphi(k)| |\widehat\mu(k)|
	\leq \left(\sum_{k = 1}^\infty |\widehat\varphi(k)| \right) \left(\sup_{j \geq 1} |\widehat\mu(j)| \right).
	$$
Now let $\chi$ be the indicator function of $[0, \varepsilon / 2]$
and take $\varphi$ to be the triangle pulse
	$$
	\varphi(x) = \left( \frac{2\chi}{\varepsilon} \right) * 
	\left( \frac{2\chi}{\varepsilon} \right),
	$$
where $*$ denotes convolution. Then
	$$
	|\widehat\varphi(k)| =
	\left| \frac{2\widehat\chi(k)}{\varepsilon}\right|^2 =
	\sinc^2\left(\frac{k\pi\varepsilon}{2} \right) \leq
	\min\left( 1, \, \frac{4}{k^2 \pi^2 \varepsilon^2} \right),
	$$
so that
	$$
	\sum_{k = 1}^\infty |\widehat\varphi(k)| \leq
	\left\lceil \frac{2}{\pi\varepsilon} \right\rceil +
	\frac{4}{\pi^2 \varepsilon^2}
	\sum_{\ceil{\frac{2}{\pi\varepsilon}} + 1}^\infty \frac{1}{k^2}
	\leq
	\frac{2}{\pi\varepsilon} + 1 +
	\frac{4}{\pi^2 \varepsilon^2} 
	\int_{\frac{2}{\pi\varepsilon}}^\infty \frac{1}{x^2} \intd{x}
	= \frac{4 + \pi\varepsilon}{\pi\varepsilon}.
	$$
It follows that
	$$
	\sup_{j \geq 1} |\widehat\mu(j)| \geq \frac{1}{2} \cdot \frac{\pi\varepsilon}{4 + \pi\varepsilon}.
	$$
\end{proof}
\end{lemma}

\begin{example}	\label{theexample}
Let $s \in (\sqrt{3} - 1, 1)$ and choose $b$ such that
	$$
	\frac{1 - s}{s} < b < \frac{s}{2}
	$$
(this is possible since $s > \sqrt{3} - 1$). Let $(l_k)_{k = 1}^\infty$ be
a sequence of natural numbers such that
	$$
	\lim_{k \to \infty} \frac{l_{k + 1}}{l_k} = \infty,
	$$
and set $m_k = \ceil{bl_k}$. Whenever $x \in [0, 1]$ is not a dyadic rational, and
thus has a unique binary decimal expansion $x = 0.x_1x_2\ldots$, let
	$$
	f(x) = \sup \left( \left\{k; \, x_{l_k + 1} \ldots x_{l_k + m_k} = 0^{m_k} \right\}
	\cup \{ 0 \} \right).
	$$
Then for each $k$
	$$
	\leb \{ x; \, f(x) = \infty \} \leq
	\leb \left( \bigcup_{j = k}^\infty \{x; \,  x_{l_j + 1} \ldots x_{l_j + m_j} = 0^{m_j}\} \right)
	\leq
	\sum_{j = k}^\infty 2^{-m_j},
	$$
where \leb denotes Lebesgue measure. The sum converges since $(m_k)$ is eventually strictly
increasing, and thus $f(x)$ is finite for \leb-a.e.~$x \in [0, 1]$. Let
	$$
	A = \{ x \in [0, 1]; \, f(x) \text{ is even} \}, \qquad
	B = \{ x \in [0, 1]; \, f(x) \text{ is odd} \};
	$$
then $\dimf (A \cup B ) = 1$ since $\leb(A \cup B) = 1$.

To see that $\dimf A \leq s$, take any $\mu \in \probs(A)$ and
define for odd $k$
	\begin{align*}
	A_k &= \{ x \in A; \, x_{l_k + 1} \ldots x_{l_k + m_k} = 0^{m_k} \} \\
	A_k^j &= \{ x \in A_k; \, f(x) = j \}.
	\end{align*}
If $x \in A_k$ then $f(x) \geq k$, so there must be some even
number $j \geq k + 1$ such that $f(x) = j$. Hence for each $k$
	$$
	A_k = \bigcup_{\substack{ j \geq k + 1 \\ j \text{ even}}} A_k^j.
	$$
Let
	$$
	\alpha_k = \mu(A_k), \qquad
	\alpha_k^j = \mu(A_k^j),
	$$
and let $P$ be the set of natural numbers $k$ such that
	$$
	\alpha_k^j \leq \frac{2^{-(m_k + j - k)}}{6} \qquad
	\text{for all even } j \geq k + 1.
	$$

Suppose first that the set $P$ is infinite. Let
$\mu_k = \restr{\mu}{A \setminus A_k}$ and let $\nu_k$ be the
image of $\mu_k$ under the map $x \mapsto 2^{l_k} x \pmod{1}$.
The measure $\nu_k$ is concentrated on $[2^{-m_k}, 1]$ and has total
mass $1 - \alpha_k$, so by Lemma~\ref{effinfsuplemma} there is some
$r_k \geq 1$ such that
	$$
	\widehat\mu_k\left(2^{l_k} r_k\right) =
	\widehat\nu_k(r_k) \geq (1 - \alpha_k) \frac{2^{-m_k}}{5}.
	$$
For $k \in P$,
	$$
	\alpha_k =
	\sum_{\substack{ j \geq k + 1 \\ j \text{ even}}} \alpha_k^j \leq \frac{2^{-m_k}}{6},
	$$
and hence
	\begin{align*}
	\left( 2^{l_k} r_k \right)^{s / 2}
	\left| \widehat\mu\left( 2^{l_k} r_k \right) \right| &\geq
	2^{sl_k / 2} \left(
	|\widehat{\mu}_k(2^{l_k} r_k)| - \alpha_k
	\right) \\
	&\geq
	2^{s l_k / 2 - m_k} \left(
	\frac{1 - \alpha_k}{5} - \frac{1}{6}
	\right),
	\end{align*}
where the exponent is positive for large $k$ since $b < s / 2$.
Thus if $(k_i)_{i = 1}^\infty$ is an enumeration of $P$ then
	$$
	\limsup_{|\xi| \to \infty} |\xi|^{s / 2} |\widehat\mu(\xi)| \geq
	\lim_{i \to \infty} \left( 2^{l_{k_i}} r_{k_i} \right)^{s / 2}
	\left| \widehat\mu\left( 2^{l_{k_i}} r_{k_i} \right) \right|
	= \infty,
	$$
so $\dimf \mu \leq s$.

Suppose on the other hand that $P$ is finite. If $k$ is odd
and $j \geq k + 1$ is even, then
	$$
	A_k^j \subset \{ x \in [0, 1]; \, x_{l_k + 1} \ldots x_{l_k + m_k} = 0^{m_k}
	\text{ and } x_{l_j + 1} \ldots x_{l_j + m_j} = 0^{m_j} \}
	\subset
	\bigcup_p I_p,
	$$
where $\{ I_p \}$ are $2^{l_j - m_k}$ intervals, each of length $2^{-(l_j + m_j)}$.
Thus
	\begin{align*}
	I_s(\mu) &\geq \sum_{p = 1}^{2^{l_j - m_k}}
	\iint_{I_p \times I_p} \left| x - y \right|^{-s} \intd{\mu}(x) \intd{\mu}(y)\\
	&\geq
	2^{s(l_j + m_j)} \sum_{p = 1}^{2^{l_j - m_k}} \mu(I_p)^2
	\geq
	2^{s(l_j + m_j)} \frac{(\alpha_k^j)^2}{2^{l_j - m_k}},
	\end{align*}
where the inequality $\| \cdot \|_2^2 \geq d^{-1} \| \cdot \|_1^2$ for norms in
$\rea^d$ is used in the last step.
If $k \notin P$ then $j = j(k)$ can be chosen such that the last expression is
greater than or equal to
	$$
	\frac{1}{36} \cdot 2 \:\widehat{\phantom{;}}\, \big(( s(1 + b) - 1 ) l_j - 2j - m_k\big).
	$$
Because $b > (1 - s) / s$ and $l_j$ grows exponentially with $j$, there is some
$\varepsilon > 0$ such that the exponent is at least
	$$
	\varepsilon l_j - m_k \geq \varepsilon l_{k + 1} - \lceil bl_k \rceil
	$$
whenever $k$ (and hence $j$) is large enough. Since $P$ contains arbitrarily large odd $k$
it follows that $I_s(\mu) = \infty$, and thus $\dimf \mu \leq s$ in this case too.

This shows that $\dimf A \leq s$, and similar computations show that $\dimf B \leq s$
as well. Thus the Fourier dimension is not finitely stable.

Consider now the $F_\sigma$-sets
	$$
	A' = \bigcup_{k \text{ even}} \closure{f^{-1}(k)}, \qquad
	B' = \bigcup_{k \text{ odd}} \closure{f^{-1}(k)},
	$$
which are slightly larger than $A$ and $B$. From the inclusion
	$$
	\closure{f^{-1}(k)} \subset
	\closure{\{ x; \, x_{l_k + 1} \ldots x_{l_k + m_k} = 0^{m_k}\}} \cap
	\bigcap_{j = k + 1}^\infty \closure{\{ x; \, x_{l_j + 1} \ldots x_{l_j + m_j} \neq 0^{m_j}\}},
	$$
it follows that the difference sets $A' \setminus A$ and $B' \setminus B$ only contain
dyadic rationals, and in particular they are countable. Thus $A'$ and $B'$ have the same
Fourier dimensions as $A$ and $B$ respectively (using that any measure that gives positive
mass to a countable set has Fourier dimension $0$), and in particular
	$$
	\dimf\left(A' \cup B'\right) > \max\left(\dimf A', \, \dimf B' \right).
	$$
On the other hand,
Proposition~5 in \cite{EPS} implies that the compact Fourier dimension is finitely
stable on $F_\sigma$-sets, and thus
	$$
	\dimfc\left(A' \cup B'\right) = \max\left(\dimfc A', \, \dimfc B'\right).
	$$
\end{example}

\bibliographystyle{plain}
\bibliography{references}

\begin{thebibliography}{1}

\bibitem{EPS}
Fredrik Ekstr{\"o}m, Tomas Persson, and J{\"o}rg Schmeling.
\newblock On the {F}ourier dimension and a modification.
\newblock Preprint: arXiv:1406.1480v3 [math.FA]. To appear in {J}ournal of
  {F}ractal {G}eometry.

\bibitem{mattila}
Pertti Mattila.
\newblock {\em Geometry of {S}ets and {M}easures in {E}uclidean {S}paces ---
  {F}ractals and {R}ectifiability}.
\newblock Cambridge {U}niversity {P}ress, 1995.

\end{thebibliography}
\end{document}